\documentclass[11pt]{amsart}
\usepackage{amssymb,amsmath}
\usepackage{latexsym}

\def\dfrac#1#2{{\displaystyle\frac{#1}{#2}}}

\newtheorem{thm}{Theorem}
\newtheorem{lem}{Lemma}

\begin{document}

\title{Growth sequences for  
circle diffeomorphisms}
\author{Nobuya Watanabe}
\address{ Department of Mathematics, School of Commerce, Waseda University,
 Shinjuku, Tokyo 169-8050, Japan.}
\email{nobu@waseda.jp}
\begin{abstract}
We obtain results on the growth sequences 
of the differential   
for iterations of circle diffeomorphisms without periodic points. 
\end{abstract}

\maketitle
\section{Introduction and statement of results}
\noindent
Let $f:S^1 \rightarrow S^1$ be a
 $C^{1}$-diffeomorphism where $S^1={\Bbb R}/{\Bbb Z}$. 
We define the {\it growth sequence} for $f$ by
\[
\Gamma_n(f)= 
\max \{ \lVert Df^n\rVert, 
\lVert Df^{-n}\rVert \},\ \ \ n \in \Bbb{N},\]
where $f^n$ is the $n$-th iteration of $f$ and 
$\lVert Df^n\rVert = {\displaystyle \max_{x \in S^1}}|Df^n(x)|$.

If $f$ has periodic points, then the study of growth sequences 
reduces to the case of interval 
diffeomorphisms which was studied in \cite{B},\cite{PS},\cite{W}.  

If $f$ has no periodic points, then by the 
theorem of Gottschalk-Hedlund $\Gamma_n(f)$ is bounded if and only if $f$ 
is $C^1$-conjugate to a rotation. Notice that if $\Gamma_n(f)$ 
is bounded then $f$ is minimal. So it is natural to ask how rapidly could 
the sequence $\Gamma_n(f)$ grow if it is unbounded.
  
In this paper we give an answer to this question:

\begin{thm}
Let $f:S^1 \rightarrow S^1$ be a $C^2$-diffeomorphism without 
periodic points.
Then 
\[\lim_{n\rightarrow \infty}\frac{\Gamma_n(f)}{n^2}=0.\]
\end{thm}

\begin{thm}
For any increasing unbounded sequence of positive real numbers 
$\theta_n =o(n^2)$ 
as $n \rightarrow \infty$ 
and any $\varepsilon >0$ there exists an analytic diffeomorphism 
$f:S^1 \rightarrow S^1$ 
without periodic points such that 
\[1- \varepsilon \le \limsup_{n\rightarrow \infty}
\frac{\Gamma_n(f)}{\theta_n} \le 1.\]
\end{thm}

\section{Preliminaries}

\noindent
Given an orientation preserving homeomorphism $f:S^1\rightarrow S^1$, 
its {\it rotation number} is defined by 
\[
\rho (f) = \lim_{n \rightarrow \infty}\frac{\tilde{f}^n(x)-x}{n} 
\mod {\Bbb Z}
\]
where $\tilde{f}$ denotes a lift of $f$ to $\Bbb{R}$. The limit exists 
and is independent on $x \in \Bbb{R}$ and a lift $\tilde{f}$. 

Put $\alpha = \rho(f)$. 
Let $R_{\alpha}$ be the rigid rotation by $\alpha$
\[
R_{\alpha}(x)=x+\alpha \mod {\Bbb Z}.
\]

For the basic properties of circle homeomorphisms and 
the combinatorics of orbits of the rotation of the circle, 
general references are 
\cite{MS} chapter I and \cite{KH} chapter 11, 12.

By Poincar$\acute{{\rm e}}$ the order structure of orbits of $f$ and 
$R_{\alpha}$ on $S^1$ 
are almost same. In particular if $\rho (f) = \frac{p}{q} \in \Bbb{Q}/\Bbb{Z}$ 
then $f$ has periodic points of period $q$ and every periodic orbits of $f$ 
have the same order as orbits of 
$R_{\frac{p}{q}}$ on $S^1$. $\rho (f) \in 
(\Bbb{R}\setminus \Bbb{Q})/\Bbb{Z}$ 
if and only if $f$ has no periodic points, in this case, if $f$ is of class 
$C^2$ then by the well known 
theorem of Denjoy $f$ is topologically conjugate to $R_{\alpha}$.

Suppose $\alpha \in (\Bbb{R}\setminus \Bbb{Q})/\Bbb{Z}$. Let
\[
\alpha = [a_1,a_2,a_3,\ldots]=\dfrac{1}{a_1+
\dfrac{1}{a_2+\dfrac{1}{a_3+{}_{\ddots}}}}\ \ \ ,  
a_{i} \ge 1, a_i \in \Bbb{N}
\]
be the continued fraction expansion of $\alpha$, and 
\[
\frac{p_n}{q_n}=[a_1,a_2,\ldots,a_n]
\]
be its $n$-th convergent. 
Then $p_n$ and $q_n$ satisfy 
\[
p_{n+1}=a_{n+1}p_n+p_{n-1},~~p_0=0,~p_1=1,
\]\vspace{-7mm}
\[
q_{n+1}=a_{n+1}q_n+q_{n-1},~~q_0=1,~q_1=a_1, 
\]
\[
\frac{p_0}{q_0} < \frac{p_2}{q_2} < \frac{p_4}{q_4} < \cdots < \alpha < \cdots 
< \frac{p_5}{q_5} < \frac{p_3}{q_3} < \frac{p_1}{q_1}.
\]

The sequence of rational numbers $\{ \frac{p_n}{q_n}\}$ is the best rational 
approximation of $\alpha$. 
This can be expressed using the dynamics of 
$R_{\alpha}$ as follows. 
$R^{q_n}_{\alpha}(0) \in [0,R^{-q_{n-1}}_{\alpha}(0)]$, and 
if $k >q_{n-1}$, 
$R^{k}_{\alpha}(0) \in [R^{q_{n-1}}_{\alpha}(0),R^{-q_{n-1}}_{\alpha}(0)]$ 
then $k \ge q_{n}$. 
Note that for $0 \le k \le a_{n+1}, 
R_{\alpha}^{kq_n}(0) \in [0,R^{-q_{n-1}}_{\alpha}(0)]$, and 
$R_{\alpha}^{(a_{n+1}+1)q_n}(0) \notin [0,R_{\alpha}^{-q_{n-1}}(0)]$.

For $\alpha \in (\Bbb{R}\setminus \Bbb{Q})/\Bbb{Z}$ 
the continued fraction expansion is unique. 
On the other hand for $\beta \in \Bbb{Q}/\Bbb{Z}$ 
expressions by continued fractions are not unique, 
$\beta = [b_1,b_2,\ldots,b_n+1] = [b_1,b_2,\ldots,b_n,1]$. 
 
For $\alpha = [a_1,a_2,\ldots]$ and $i, j \in \Bbb{N}, 
1 \le i \le j$ we denote 
$\alpha|[i,j] = [a_i,a_{i+1},\ldots,a_j]$.
In case we emphasize $\alpha$ we denote $a_{i}(\alpha), p_i(\alpha),
q_i(\alpha)$.

For $x \in S^1$, $I_n(x)$ denotes the smaller interval with endpoints $x$ and 
$f^{q_n}(x)$  and for an interval $J \subset S^1$, $|J|$ the length of $J$. 

The following is well known. See \cite{MS} chapter I section 2a. 

\begin{lem} {\rm (Denjoy)} 
 Let $f$ be a $C^1$-diffeomorphism of $S^1$ without periodic points and 
$\log Df:S^1\rightarrow {\Bbb R}$ has bounded variation. 
Then there exists
 a positive constant $C_1 = C_1(f)$  satisfying the following properties. 

$(1)$~ For any $0\le l \le q_{n+1}$
 and for every $x_1,x_2 \in I_n(x)$ 
\[
\frac{1}{C_1}\le \frac{Df^l(x_1)}{Df^l(x_2)}\le C_1. 
\]

$(2)$~{\rm (Denjoy inequality)} 
 For every $n \in {\Bbb N}$,
\[
\frac{1}{C_1}\le \lVert Df^{q_n} \rVert \le C_1. 
\]
\end{lem}

As stated in section 1, the growth sequences play a significant role in 
the problem of the smooth linearization of circle diffeomorphisms, 
where the arithmetic property of rotation numbers and the regularity of 
diffeomorphisms are important.    
This problem has a rich history, 
see e.g. \cite{A}, \cite{H}, \cite{Y}, \cite{KS}, \cite{St}, \cite{KO}.
 
In this paper, particularly we need the following 
improvement of Denjoy inequality which 
is due to Katznelson and Ornstein. The statement of Lemma 2 
is obtained by merging 
results in \cite{KO}, for $(1)$, (1.16), 
lemma 3.2 (3.6) and proposition 3.3 (a), for $(2)$,  
theorem 3.7.

\begin{lem} Let  $f$  be a
$C^2$-diffeomorphism of $S^1$ without periodic points. 
 Set 
\[
E_n = \max \{ \lVert \log Df^{q_n}\rVert,~ 
\max_{x\in S^1} \{|D\log Df^{q_n}(x)||I_{n-1}(x)|\} \}.
\]
Then the following hold.

$(1)$~~~$\lim_{n \rightarrow \infty} E_n = 0.$

$(2)$ If $f$ is of class $C^{2+\delta}, \delta > 0$ then 
there exist $C>0$ and $0 < \lambda <1$ such that 
$\lVert \log Df^{q_n}\rVert \le C\lambda^n$ 
for any $n \in \Bbb{N}$.
\end{lem}

The conclusion of Lemma 2 (2) plus some arithmetic condition of $\rho (f)$ 
are sufficient to provide the $C^{1}$-linearization of $f$. We need the 
following which is a special case of the main theorem in \cite{KO}. 
For $C^{3+\delta}$-diffeomorphisms it is originally due to Herman \cite{H}.
\vspace{3mm}

\noindent
{\bf Corollary of Lemma 2 (2).}\ \ {\it If} $f$ {\it is of class} 
$C^{2+\delta}$ {\it and the rotation number} 
$\alpha = \rho(f)$ {\it is of bounded type i.e.} 
$a_i(\alpha)$ {\it is uniformly bounded then} 
$\lVert Df^n \rVert$ {\it is uniformly bounded.} 

\section{Proof of Theorem 1}
\noindent
Let $f:S^1 \rightarrow S^1$ be a $C^2$-diffeomorphism without 
periodic points with the rotation number $\rho(f) = [a_{1},a_{2},\ldots]$ and 
its convergents $\{\frac{p_{n}}{q_{n}}\}$.

The following crucial and fundamental lemma is due to Polterovich and Sodin 
(\cite{PS} lemma 2.3).

\begin{lem} {\rm (Growth lemma)} Let $\lbrace A(k) \rbrace_{k\ge 0}$
 be a sequence of real numbers such that for each $k \ge 1$
\[
2A(k)-A(k-1)-A(k+1)\le C\exp (-A(k)),\quad C > 0,
\]
 and $A(0)=0$.  Then either for each $k \ge 0$
\[
A(k) \le 2\log \left ( k\sqrt{\frac{C}{2}}+1 \right ),
~ or
~~\liminf_{k\rightarrow \infty}\frac{A(k)}{k}>0.
\]
\end{lem}

\begin{lem}  For $0 \le k \le a_{n+1}+1$ we set 
$A_n(k)=\log \lVert Df^{kq_n}\rVert$. Then there exists a positive 
constant $C=C(f)$ 
independent with $n$ such that for $1\le k \le a_{n+1}$,
\[
2A_n(k)-A_n(k-1)-A_n(k+1)\le CE_n \exp(-A_n(k)).
\]

\end{lem}

\begin{proof}~~Let $A_n(k) = \log Df^{kq_n}(x_0)$ and $x_i=f^{iq_n}(x_0)$.
 Then we have,
\[
2A_n(k)-A_n(k-1)-A_n(k+1) 
\]
\[
\le 
2\log Df^{kq_n}(x_0)-\log Df^{(k-1)q_n}(x_1)
-\log Df^{(k+1)q_n}(x_{-1})
\]
\[
\le |\log Df^{q_n}(x_0)-\log Df^{q_n}(x_{-1})|
= |D\log Df^{q_n}(y_0)||I_n(x_{k-1})|\frac{|I_n(x_{-1})|}{|I_n(x_{k-1})|} ,
\]
where $y_0\in I_n(x_{-1})$.

Notice that the intervals $I_n(x_{-1}), I_n(x_{0}), I_n(x_{1}),\ldots, 
I_n(x_{a_{n+1}-1})$ are adjacent in this order 
and ${\displaystyle \cup_{i=0}^{a_{n+1}-1}}
I_n(x_i) \subset 
I_{n-1}(f^{-q_{n-1}}(x_0))$. 
Since $y_0 \in I_n(x_{-1})$, we have for $1\le k \le a_{n+1}-1$, 
$I_n(x_{k-1})\subset I_{n-1}(f^{-q_{n-1}}(y_0))$. So by Denjoy inequality 
(Lemma 1 (2)) we have 
\[
|I_n(x_{k-1})|\le C_1^2|I_{n-1}(y_0)|,
\]
and using lemma 1 (1) we have
\[
\frac{|I_n(x_{-1})|}{|I_n(x_{k-1})|} \le C_1\frac{1}{Df^{kq_n}(x_0)}.
\]
Hence we have
\[
2A_n(k)-A_n(k-1)-A_n(k+1) 
\]
\[
\le C_1^3|D\log Df^{q_n}(y_0)||I_{n-1}(y_0)| \frac{1}{Df^{kq_n}(x_0)} 
\le C_1^3E_n\exp(-A_n(k)).
\]

\end{proof}

We extend $A_n(k)$ for $k \ge a_{n+1}+2$ by $A_n(k)=A_n(a_{n+1}+1)$. 
Then by Lemma 1 (2) and the definition of $E_{n}$ we have
\[
2A_n(a_{n+1}+1)-A_n(a_{n+1})-A_n(a_{n+1}+2)
\]
\[\le \log Df^{(a_{n+1}+1)q_n}(x_0)-\log Df^{a_{n+1}q_n}(x_0)
\le \lVert \log Df^{q_n}\rVert
\]
\[
\le E_n \exp (-A_n(a_{n+1}+1))\lVert Df^{(a_{n+1}+1)q_n}\rVert
\]
\[
\le E_n \exp (-A_n(a_{n+1}+1))\lVert Df^{q_{n+1}}\rVert 
\lVert Df^{q_{n}}\rVert \lVert Df^{-q_{n-1}}\rVert
\]
\[
\le C_1^3E_n \exp (-A_n(a_{n+1}+1)).
\]

For $k \ge a_{n+1}+2$, $2A_n(k)-A_n(k-1)-A_n(k+1)=0$.

Then since $A_n(k)$ satisfy the condition of Lemma 3 with the constant 
$C=C_1^3$ and  obviously $ \lim_{k \rightarrow \infty}\frac{A_n(k)}{k}=0$,  
we have 
\[
\lVert Df^{kq_n}\rVert \le \left (\sqrt{\frac{CE_n}{2}}k+1 \right) ^2,
~~0\le k \le a_{n+1}.\]

For $q_n\le l < q_{n+1}$, we define $0\le k_{i+1}\le a_{i+1}, 
(i=0,1,\ldots,n)$ 
inductively by
\[
r_{n+1}=l,~~ r_{i+1}=k_{i+1}q_i+r_i,~~ 0\le r_i < q_i.
\]

Then, using ~$\frac{q_{i+1}}{q_i} \ge a_{i+1} \ge k_{i+1}$,
\[
\frac{\lVert Df^l\rVert}{l^2}\le 
\frac{\prod^n_{i=0}\lVert Df^{k_{i+1}q_i}\rVert}{(k_{n+1}q_n)^2}
\le 
\frac{\prod^n_{i=0}\left (\sqrt{\frac{CE_i}{2}}
 k_{i+1}+1\right )^2}
{\left (k_{n+1}\prod^{n-1}_{i=0}\frac{q_{i+1}}{q_i}\right )^2}
\]
\[
\le \left ( \sqrt{\frac{CE_n}{2}}+1\right ) ^2
\prod^{n-1}_{i=0}\left( \sqrt{\frac{CE_i}{2}}
+\frac{q_i}{q_{i+1}}\right)^2.
\]

Since $\frac{q_i}{q_{i+2}}<\frac{1}{2}$, for sufficiently small $E_i$ and
 $E_{i+1}$ 
\[
\left( \sqrt{\frac{CE_i}{2}}+\frac{q_i}{q_{i+1}}\right) 
\left( \sqrt{\frac{CE_{i+1}}{2}}+\frac{q_{i+1}}{q_{i+2}}\right) 
\le \frac{1}{2}.
\]
By Lemma 2 (1), $E_n \rightarrow 0$ as $n \rightarrow \infty$.
 Consequently we have
\[
\lim_{l \rightarrow \infty}\frac{\lVert Df^l\rVert}{l^2}=0.
\]

For the case $\lVert Df^{-l} \rVert, l >0$, the argument is the same.

\section{Proof of Theorem 2}
\noindent
Let $\{\theta_n\}_{n\ge 1}$ be any increasing unbounded sequence of positive 
real numbers such that 
$\theta_n=o(n^2)$ as $n \rightarrow \infty$.

We consider the two-parameter family of rational functions 
on the Riemann sphere 
$\hat{\Bbb{C}}= \Bbb{C}\cup \{\infty\}$,
\[
J_{a,t}: \hat{\Bbb{C}} \rightarrow \hat{\Bbb{C}}, \ \ 
J_{a,t}(z)= \exp (2\pi i t)z^2\frac{z+a}{az+1} 
\]
where $a \in \Bbb{R}, a > 3$ and $t \in \Bbb{R}/\Bbb{Z}$. 

For each $a, t$ the map $J_{a,t}$ makes invariant the unit circle 
$\partial \Bbb{D} = \{z\in \Bbb{C}; \lvert z \rvert =1\}, 
J_{a,t}(\partial \Bbb{D})= \partial \Bbb{D}$, moreover the restriction of 
$J_{a,t}$ to $\partial \Bbb{D}$ is an orientation preserving diffeomorphism.  
The set of critical points of $J_{a,t}$ consists of four elements containing 
$0$ and $\infty$ which are fixed by $J_{a,t}$. Notice that if 
$a \rightarrow \infty$ then on a compact tubular neighbourhood of 
the unit circle in $\Bbb{C}\setminus \{0\}$ $J_{a.t}$ uniformly converges to the rotation $z 
\mapsto \exp (2\pi i t)z$. 

Put $\psi : \Bbb{R}/\Bbb{Z} \rightarrow \partial\Bbb{D}, 
\psi (x)= \exp (2\pi i x)$. Conjugating $J_{a,t}|\partial \Bbb{D}$ by 
$\psi$ we obtain the family of analytic circle diffeomorphisms 
$\{f_{a,t}\}$,
\[
f_{a,t}: \Bbb{R}/\Bbb{Z} \rightarrow \Bbb{R}/\Bbb{Z},~~
f_{a,t}(x)=\psi^{-1}\circ J_{a,t}\circ \psi (x) 
=f_{a,0}(x)+t \mod{\Bbb{Z}}.
\] 
Temporarily we fix $a > 3$ and abbreviate as $f_{a,t}= f_{t}$.

The following properties of this family are standard. 
See e.g. \cite{MS} chapter I, section 4, where Arnold family 
$x \mapsto x+ a\sin (2\pi x) +t$ is mainly 
dealt with but the argument is valid for our family. 
Also see \cite{KH} chapter 11, section 1.

The map $F: S^1 \rightarrow S^1 , t \mapsto \rho (f_t)$ is continuous 
and monotone increasing. We set
\[
K=\{t\in S^1; \rho (f_t) ~~{\rm is~ irrational} \}.
\]
We denote Cl($K$) the closure of 
$K$.  $F|K$ is a one-to-one map. 
For $t \in K$ with $F(t)=\alpha$, we denote $f_{t}=\hat{f}_{\alpha}$. 
Notice that $f_{t}$ never conjugate to a rational rotation. 
Hence for $\frac{p}{q} \in \Bbb{Q}/\Bbb{Z}$, $F^{-1}(\frac{p}{q})$ 
is a closed interval, say, 
$[\frac{p}{q}_{-},\frac{p}{q}_{+}] $. 

Moreover, $F^{-1}|(\Bbb{R}\setminus \Bbb{Q})/\Bbb{Z}: 
(\Bbb{R}\setminus \Bbb{Q})/\Bbb{Z} 
\rightarrow K$ is continuous and 
\[
\lim_{\alpha \rightarrow \frac{p}{q}-0}F^{-1}|
(\Bbb{R}\setminus \Bbb{Q})/\Bbb{Z}(\alpha)
 = \frac{p}{q}_{-}, 
 ~\lim_{\alpha \rightarrow 
 \frac{p}{q}+0}F^{-1}|(\Bbb{R}\setminus \Bbb{Q})/\Bbb{Z}(\alpha)
 = \frac{p}{q}_{+}.
\] 

Note that for every $\frac{p}{q}\in \Bbb{Q}/\Bbb{Z}$ and every $x \in S^{1}$, 
there exists $t \in [\frac{p}{q}_{-},\frac{p}{q}_{+}]$ such that 
$f_{t}^{q}(x)=x$.
For $\frac{p}{q} \in \Bbb{Q}/\Bbb{Z}$, put $t_{*} = \frac{p}{q}_{-}$. 
The case $t_{*}= \frac{p}{q}_{+}$ is similar. 
Then the graph of $f^q_{t_{*}}(x)$ 
touches from below to the graph of the identity map, 
in particular, there exists $x_0 \in S^1$ such that 
\[
f^q_{t_{*}}(x_0)=x_0, ~~Df^q_{t_{*}}(x_0)=1.
\]
Then the following holds. 
\begin{lem}
$D^2f^q_{t_{*}}(x_0) \ne 0.$
\end{lem}
\begin{proof}
By contradiction, we suppose $D^2f_{t_{*}}^q(x_0)=0$. 
Then in our case $D^3f_{t_{*}}^q(x_0)=0$, otherwise $x_0$ is a topologically 
transversal fixed point of $f^q_{t_{*}}$ and persists under 
perturbation of $f_{t_{*}}$, which 
contradicts $t_{*} \in \mathrm{Cl}(K)\setminus K$. Set $z_{0} = 
\psi (x_{0}) \in \partial \Bbb{D}$. Since the order of tangency to the 
identity map  
is an invariant of $C^{\infty}$-conjugacy \cite{T},   
we have for $J_{t_{*}}=J_{a,t_{*}}$ 
\[
J_{t_{*}}^{q}(z_{0})=z_{0}, \ DJ_{t_{*}}^{q}(z_{0})=1, 
\ D^{2}J_{t_{*}}^{q}(z_{0})= 
D^{3}J_{t_{*}}^{q}(z_{0})= 0.
\]
So $z_{0}$ is a parabolic fixed point for $J_{t_{*}}^{q}$ with multiplicity 
at least four. See \cite{M} chapter 7. By the Laeu-Fatou flower theorem 
(\cite{M} th.7.2) $z_0$ has at least three basins of attraction for 
$J^q_{t_{*}}$. Let $B$ be 
 one of the immediate attracting basins of $z_0$ for $J^q_{t_{*}}$. 
Then $B$ must contain at least one critical point 
of $J_{t_{*}}^{q}$ (\cite{M} corollary 7.10). 
So each basin of the cycle $\{z_0,J_{t_{*}}(z_0),\ldots,
 J_{t_{*}}^{q-1}(z_{0})\}$  contains 
 at least one critical point of $J_{t_{*}}$. But $J_{t_{*}}$ has exactly 
four critical points and two of them are fixed points. 
We obtain a contradiction.
\end{proof}

Hence, for example, by comparing a fractional linear transformation 
(see also \cite{B} thorem 1 (A)), 
we can see that there exist $C >0$ and 
$\{x_l\}_{l\ge 1} \subset S^1 $ with 
$\lim_{l\rightarrow \infty}x_l=x_0$ such that 
\[
Df^{lq}_{t_{*}}(x_l)\ge Cl^2 , ~~
\mathrm{for~any}~l\in \Bbb{N}.
\]
Since $\theta_n=o(n^2)$ , we have

\vspace{3mm}
\noindent
{\bf Corollary of Lemma 5.}\ \ {\it For sufficiently large} $l$, {\it we have} 
$\lVert Df^{lq}_{t_{*}} \rVert > \theta_{lq}$. 
\vspace{3mm}

\noindent
{\bf Remark.}  For each $k \in {\Bbb N}$ we set
\[
U_k=\{ t\in \mathrm{Cl}(K); \mathrm{There~ exist}\ m\ge k ~
\mathrm{and}~ x\in S^1
\ \mathrm{such\ that}\  Df^m_t(x) > m\sqrt{\theta_m} \}.
\]
Obviously $U_k$ is open set in Cl($K$). By the corollary and the denseness 
of preimages of rational numbers by $F$ in Cl($K$), 
$U_k$ is dense in Cl($K$). So the following set is a residual subset 
of Cl($K$), 
\[
\{t \in \mathrm{Cl}(K); \limsup_{n\rightarrow \infty}
\frac{\Gamma_n(f_{t})}{\theta_n}=\infty\}.
\]
\vspace{1mm}

We seek a desired diffeomorphism in this family $\{ f_t \}$ by specifying its 
rotation number $\alpha_{\infty} = \rho(f_{t_{\infty}}) 
\in (\Bbb{R}\setminus \Bbb{Q})/\Bbb{Z}$.
We will define an increasing sequence of even numbers 
$0 < n_1 < n_2 < n_3< \cdots$, and 
a sequence of positive integers $A_1,A_2,A_3,\ldots$ inductively. 
The continued fraction expansion of $\alpha_{\infty}$ is the following. 
\[
\alpha_{\infty}
=[a_1(\alpha_{\infty}),a_2(\alpha_{\infty}),a_3(\alpha_{\infty}),\ldots] 
\]
\[
= [1,1,\ldots,1,A_{1},1,\ldots,1,A_{2},1,\ldots,1,A_{k},1,\ldots]
\]
where if $i=n_k$ then $a_i(\alpha_{\infty})=A_k$ 
and if $i \ne n_k$ for any $k$ then $a_i(\alpha_{\infty})=1$. 

For $m, A \ge 1, m, A \in \Bbb{N}$, we set
\[
\alpha_{m}^{A} 
= [a_1(\alpha_{m}^{A}),a_2(\alpha_{m}^{A}),a_{3}(\alpha_{m}^{A}),\ldots]
\]
\[
= [1,1,\ldots,1,A_{1},1,\ldots,1,A_{m-1},1,\ldots,1,A,1,1,1,\ldots]
\]
where $a_i(\alpha_{m}^{A})=A_k$ if 
$i= n_k \le n_{m-1}$ and $a_i(\alpha_{m}^{A})=A$ 
if $i=n_m$ and $a_i(\alpha_{m}^{A})=1$ otherwise. 

Set $\alpha_m = \alpha_{m}^{A_m}$.
Notice that $\alpha_{m}^{A}|[1,n_{m}-1] = \alpha_{\infty}|[1,n_{m}-1]$ 
and $\alpha_{m}^{A}$ is of bounded type.
Unless otherwise 
stated we use the symbols
$p_{n}, q_{n}$ as $p_{n}(\alpha_{\infty}), q_{n}(\alpha_{\infty})$.

\begin{lem}
There exist a sequence of even numbers $0 < n_1 < n_2 < n_3< \cdots$, and 
a sequence of positive integers $A_1,A_2,A_3,\ldots$ 
such that for each $m \ge 1$ the following properties hold.

$(1)$ For any $j \in \Bbb{Z}$ with $q_{n_{m}-1} \le |j| \le A_mq_{n_{m}-1}$, 
$\lVert D\hat{f}_{\alpha_{m}}^{j} \rVert < \theta_{|j|}$.

$(2)$ There exists $j_{m}\in \Bbb{Z}$  such that 
\[
q_{n_{m}-1} \le |j_{m}| \le (A_m+1)q_{n_{m}-1} ,\ 
\lVert D\hat{f}_{\alpha_{m}^{A_{m}+1}}^{j_{m}} \rVert \ge \theta_{|j_{m}|}.
\]

$(3)$ For any $t \in F^{-1}(\alpha)$ with 
$\alpha|[1,n_{m+1}-1] = \alpha_{m}|[1,n_{m+1}-1]$ and 
any $j \in \Bbb{Z}$ with $\lvert j \rvert \le q_{n_{m}}$, 
\[
\lVert Df_{t}^{j} \rVert-1 \le \lVert D\hat{f}_{\alpha_{m}}^{j}\rVert \le 
\lVert Df_{t}^{j}\rVert +1 .
\]
\end{lem}
\begin{proof}
Let $\alpha_0 = [1,1,1,\ldots] = \frac{\sqrt{5}-1}{2}$. Since $\alpha_0$ is 
of bounded type by Corollary of Lemma 2 (2) there exists $C_0 > 0$ such that 
for any $l \in \Bbb{Z}$ , 
$\lVert D\hat{f}^{l}_{\alpha_{0}}\rVert \le C_{0}$. Let $n_1$ be a 
sufficiently large even 
number 
such that if $|i| \ge q_{n_1-1}(\alpha_0)$ then $\theta_{|i|} \ge C_0$. 

Let $\beta_1 = \alpha_0|[1,n_{1}-1] = 
\frac{p_{n_1-1}(\alpha_0)}{q_{n_1-1}(\alpha_0)} = 
[1,1,\ldots1]=[1,1,\ldots1,\infty] \in \Bbb{Q}/\Bbb{Z}$. Then by Corollary of Lemma 5 
there exists $d \in \Bbb{N}$ 
such that $\lVert Df_{\beta_{1-}}^{dq_{n_1-1}}\rVert > \theta_{dq_{n_1-1}}$, 
where $F^{-1}(\beta_{1})=[\beta_{1-},\beta_{1+}]$. 
Since $\alpha_{1}^{A} \rightarrow \beta_1-0$ as $A \rightarrow \infty$, 
$F^{-1}(\alpha_{1}^{A}) \rightarrow \beta_{1-}$ as 
$A \rightarrow \infty$. So for sufficiently large $A$ we have 
$\lVert D\hat{f}_{\alpha_{1}^{A}}^{dq_{n_1-1}}\rVert > \theta_{dq_{n_1-1}}$.
Hence the following is well defined.
\[
A_{1}=\max\{ A; \mathrm{for\  any\ } j \in \Bbb{Z}\ \mathrm{with}\  
q_{n_1-1} \le |j| \le Aq_{n_1-1},\  
\lVert D\hat{f}_{\alpha_{1}^{A}}^{j}\rVert < \theta_{\lvert j\rvert}\}.
\]
Therefore there exists $j_{1} \in \Bbb{Z}$ such that 
\[
q_{n_1-1} \le |j_{1}| \le (A_{1}+1)q_{n_1-1}, \ 
\lVert D\hat{f}_{\alpha_{1}^{A_{1}+1}}^{j_{1}}\rVert \ge \theta_{\lvert j_{1}\rvert}.
\]

Suppose we have $n_{1}, n_{2},\ldots,n_{m-1}$ and $A_{1},A_{2},\ldots,A_{m-1}$ 
satisfying conditions of Lemma. Notice that $\alpha_{m-1}$ is of bounded type 
and that (3) is satisfied by only requiring that 
$n_{m}-n_{m-1}$ is sufficiently large. 
So by the exactly same procedure as above 
we choose a sufficiently large even number $n_{m}$ and set
\[
A_{m}=\max\{ A; \mathrm{for\  any\ } j \in \Bbb{Z}\ \mathrm{with}\  
q_{n_m-1} \le |j| \le Aq_{n_m-1},\  
\lVert D\hat{f}_{\alpha_{m}^{A}}^{j}\rVert < \theta_{\lvert j\rvert}\}.
\]

\end{proof}

\begin{lem}
Let $\beta_{0}, \beta_{1}, \beta_{2} \in \Bbb{Q}/\Bbb{Z}$ be 
\[
\beta_{i} = [b_{1}(\beta_{i}),b_{2}(\beta_{i}),\ldots,b_{2n}(\beta_{i})] 
= \frac{p_{2n}(\beta_{i})}{q_{2n}(\beta_{i})}, \ i=0,1,2
\]
such that $\beta_{0}|[1,2n-1]=\beta_{1}|[1,2n-1]=\beta_{2}|[1,2n-1] $ 
and for some $B \ge 1, B \in \Bbb{N},\  b_{2n}(\beta_{i})=B+i$.

Then for any $s_{1}, s_{2} \in F^{-1}((\beta_{0},\beta_{2}))$ 
and any $x \in S^1$ we have
\[
\sum_{i=1}^{q_{2n}(\beta_{2})}\lvert ( f_{s_{1}}^{i}(x), 
f_{s_{2}}^{i}(x) ) \rvert \le 7.
\]
\end{lem}
\begin{proof}
The argument of the proof is same as the \'Swi\c atek's of lemma 3 in 
\cite{Sw}. 
We recall Farey interval. A Farey interval is an interval 
$I = (\frac{p}{q},\frac{p'}{q'}), p,p',q,q'\in \Bbb{Z}, q,q' > 0$  
with $pq'-p'q=1$. 
Then the following holds. 
\vspace{2mm}

$(*)$ All rational in $I$ have the form \ 
$\displaystyle \frac{kp+lp'}{kq+lq'}, \ k, l \ge 1, k,l \in \Bbb{N}$.
\vspace{2mm}

Since $q_{2n}(\beta_{i})=(B+i)q_{2n-1}(\beta_{0})+q_{2n-2}(\beta_{0})$ 
and 
$p_{2n}(\beta_{i})=(B+i)p_{2n-1}(\beta_{0})+p_{2n-2}(\beta_{0})$ 
two intervals $(\beta_{0},\beta_{1}), (\beta_{1},\beta_{2})$ are 
Farey intervals and 
$q_{2n}(\beta_{0})<q_{2n}(\beta_{1})<q_{2n}(\beta_{2})$ and by $(*)$ 
the cardinality of the set of
rationals in  $(\beta_{0},\beta_{2})$ with denominator less than 
$2q_{2n}(\beta_{2})$ is at most six (three if $B \ge 3$).

For given $x \in S^1$ we define 
\[
t_{1}=\sup\{t \in [\beta_{0-},\beta_{0+}]; f_{t}^{q_{2n}(\beta_{0})}(x)=x\},
\]
\[
t_{2}=\inf\{t \in [\beta_{2-},\beta_{2+}]; f_{t}^{q_{2n}(\beta_{2})}(x)=x\}.
\]

We define a diffeomorphism $G: S^{1}\times [t_{1},t_{2}] \rightarrow
S^{1}\times [t_{1},t_{2}]$ by $G(y, t)=(f_{t}(y),t)$. 
Then we have 
\[
DG^{i}(y,t)=
\left(
\begin{array}{@{\,}cc@{\,}}
Df_{t}^{i}(y)&\frac{d}{dt}(f_{t}^{i}(y))\\
0&1
\end{array}
\right)
=
\left(
\begin{array}{@{\,}cc@{\,}}
Df_{t}^{i}(y)&1+\sum_{k=1}^{i-1}Df_{t}^{i-k}(f_{t}^{k}(y))\\
0&1
\end{array}
\right).
\]
So $G$ monotonically twists $S^1$-direction to the right. 
More precisely, let $\tilde{G}: \Bbb{R}\times [t_{1},t_{2}] \rightarrow 
\Bbb{R}\times [t_{1},t_{2}], \tilde{G}(\tilde{y},t)= 
(\tilde{f_{t}}(\tilde{y}),t)$ be a lift of $G$, then for any $i \ge 1$ 
the slope of the image 
of a vertical segment $\{ \tilde{y}\} \times 
[t_{1},t_{2}]$ by $\tilde{G}^{i}$ 
is everywhere positive finite. 
Let $P : S^{1}\times [t_{1},t_{2}] \rightarrow S^1$ be 
the projection on the first coordinate. 

By contradiction we assume $\sum_{i=1}^{q_{2n}(\beta_{2})}
\lvert ( f_{s_{1}}^{i}(x), 
f_{s_{2}}^{i}(x) ) \rvert > 7$. 
We consider the interval $\gamma = \{x\}\times [t_{1},t_{2}]$ 
and its images by $G^{i}$. 
Since $[s_{1},s_{2}] \subset (t_{1},t_{2})$, 
intervals $P(G^{i}(\gamma)), 1 \le i \le q_{2n}(\beta_{2})$ overlap 
somewhere with multiplicity at least eight. 
Then, by the twist condition of $G$ 
there exist distinct natural numbers $i_{k}$, ($0 \le k \le 7, k \in \Bbb{Z}$) 
with $1 \le i_{k} \le q_{2n}(\beta_{2})$ such that  
for each $k$ ($1 \le k \le 7$),
\[
(\{f_{t_{2}}^{i_{0}}(x)\} \times [t_{1},t_{2}]) 
\cap G^{i_{k}}(\gamma)\ne \emptyset.
\]
Moreover, using the preservation of order by
$\tilde{f_{t}}: \Bbb{R}\times \{t\} \rightarrow \Bbb{R}\times \{t\}$ 
and the twist condition of $G$,   
we can see that for any $j \ge 0$, 
\[
(\{f_{t_{2}}^{i_{0}+j}(x)\} \times [t_{1},t_{2}]) 
\cap G^{i_{k}+j}(\gamma)\ne \emptyset.
\]
In particular for $j=q_{2n}(\beta_{2})-i_{0}$ by the definition of $t_{2}$ 
we have
\[
\gamma \cap G^{i_{k}+q_{2n}(\beta_{2})-i_{0}}(\gamma)\ne \emptyset.
\]
This imply that there exists a parameter value 
$u_{k} \in (t_{1},t_{2})$ such that 
$f_{u_{k}}^{q_{2n}(\beta_{2})+i_{k}-i_{0}}(x)=x$. 
For each $k$ ($1 \le k \le 7$) the denominator of 
$\rho(f_{u_{k}})$ which divides $q_{2n}(\beta_{2})+i_{k}-i_{0}$ 
is less than $2q_{2n}(\beta_{2})$. This is a contradiction.
\end{proof}

\noindent
{\it Proof of Theorem 2.}

\noindent
$\star$ {\bf Lower bound.} 
Let $j_{m} \in \Bbb{Z}$ be in Lemma 6 (2). 
Then $\lvert j_{m}\rvert \le (A_{m}+1)q_{n_{m-1}} < 
q_{n_{m}}(\alpha_{m}^{A_{m}+2})$. We assume $j_{m} > 0$. 
Then since three rational numbers 
\[
\alpha_{m}^{A_{m}}|[1,n_{m}], \alpha_{m}^{A_{m}+1}|[1,n_{m}], 
\alpha_{m}^{A_{m}+2}|[1,n_{m}]
\]
satisfy the condition of Lemma 7 and 
\[
\alpha_{\infty} \in 
(\alpha_{m}^{A_{m}}|[1,n_{m}], \alpha_{m}^{A_{m}+1}|[1,n_{m}]),
\]
\[  
\alpha_{m}^{A_{m}+1} \in
(\alpha_{m}^{A_{m}+1}|[1,n_{m}], \alpha_{m}^{A_{m}+2}|[1,n_{m}]),
\]
we have for any $x \in S^{1}$
\[
\lvert \log D\hat{f}_{\alpha_{\infty}}^{j_{m}}(x) - 
\log D\hat{f}_{\alpha_{m}^{A_{m}+1}}^{j_{m}}(x) \rvert
\]
\[ 
= 
\left \lvert \sum_{i=1}^{j_{m}-1}
\log Df_{0}(\hat{f}_{\alpha_{\infty}}^{i}(x)) - 
\sum_{i=1}^{j_{m}-1}\log Df_{0}(\hat{f}_{\alpha_{m}^{A_{m}+1}}^{i}(x)) \right \rvert  
\]
\[
\le
\rVert D\log Df_{0} \rVert 
\sum_{i=1}^{j_{m}-1}\
\lvert ( \hat{f}_{\alpha_{\infty}}^{i}(x), \hat{f}_{\alpha_{m}^{A_{m}+1}}^{i}(x) ) \rvert 
\le 
7\lVert D\log Df_{0} \rVert .
\]
Since there exists $x_{*} \in S^1$ such that 
$\lvert D\hat{f}_{\alpha_{m}^{A_{m}+1}}^{j_{m}}(x_{*})\rvert \ge \theta_{j_{m}}$ 
we have 
\[
\frac{\lVert D\hat{f}_{\alpha_{\infty}}^{j_{m}}\rVert}{\theta_{j_{m}}} 
\ge 
\frac{\lvert D\hat{f}_{\alpha_{\infty}}^{j_{m}}(x_{*})\rvert}
{\lvert D\hat{f}_{\alpha_{m}^{A_{m}+1}}^{j_{m}}(x_{*})\rvert } 
\ge 
\exp (-7\lVert D\log Df_{0} \rVert).
\]

For the case $j_{m} < 0$, using the chain rule 
$D\hat{f}_{\alpha}^{j_{m}}(x) = 
(D\hat{f}_{\alpha}^{-j_{m}}
(\hat{f}_{\alpha}^{j_{m}}(x)))^{-1}$ 
we can obtain the same estimates . 

As stated above by making the parameter $a$ sufficiently large 
we can assume that 
$\lVert D\log Df_{0} \rVert = \lVert D\log Df_{a,0} \rVert$ 
is smaller than any given positive value. 
\vspace{2mm}

\noindent
$\star$ {\bf Upper bound.} 
Let $l \in \Bbb{Z}$ with $q_{n} \le l < q_{n+1}$. The case 
$q_{n} \le -l < q_{n+1}$ is similar. 
Let $n_{m} =\max \{n_{i}; n_{i}\le n\}$. As in the proof of Theorem 1 
we expand $l$ as follows, 
\[
l = k_{n+1}q_{n} + \cdots + k_{n_{m}+1}q_{n_{m}}+cq_{n_{m}-1}+r, 
\]
where $0 \le k_{i} \le a_{i}(\alpha_{\infty})=1$ 
($n_{m}+1 \le i \le n+1$) and we choose 
$c \in \{-1,0,1\}$ so that 
$q_{n_{m}-1} \le r \le A_{m}q_{n_{m}-1}$. 

By Lemma 2 (2) and Lemma 6 (1), (3) we have 
\[
\lVert D\hat{f}_{\alpha_{\infty}}^{l}\rVert \le 
\lVert D\hat{f}_{\alpha_{\infty}}^{q_{n}}\rVert \cdots 
\lVert D\hat{f}_{\alpha_{\infty}}^{cq_{n_{m}-1}}\rVert
\lVert D\hat{f}_{\alpha_{\infty}}^{r}\rVert 
\]
\[
\le 
\exp(C\sum_{i=n_{m}-1}^{n}\lambda^{i})(1+\lVert D\hat{f}_{\alpha_{m}}^{r}\rVert) 
\le
\exp(C\sum_{i=n_{m}-1}^{n}\lambda^{i})(1+\theta_{r}). 
\]
Therefore we have 
\[
\limsup_{l \rightarrow \infty}
\frac{\lVert D\hat{f}_{\alpha_{\infty}}^{l}\rVert}{\theta_{l}} 
\le
\limsup_{l \rightarrow \infty}
\frac{\exp(C\sum_{i=n_{m}-1}^{n}\lambda^{i})(1+\theta_{r})}{\theta_{l}} 
\le 1.
\]


\vspace{4mm}

\begin{thebibliography}{[AA]}

\bibitem[A]{A}V. I. Arnold, Small denominators I, on the mapping of 
a circle into itself, Izv. Akad. Nauk. serie Math. 25 (1) (1961), 21-86. 
Translation Amer. Math. Soc. 2nd series, 46,213-284.

\bibitem[B]{B}A. Borichev, Distortion growth for iterations of diffeomorphisms 
of the interval, Geom.Funct.Anal. 14(2004), no.5, 941-964.

\bibitem[H]{H}M. R. Herman, Sur la conjugation diff$\acute{{\rm e}}$rentiable 
des diff$\acute{{\rm e}}$omorphismes du cercle $\grave{{\rm a}}$ des rotations, Publ. Math. I.H.E.S. 49 (1979), 5-234.  

\bibitem[KH]{KH}A. Katok, B. Hasselblatt, Introduction to the Modern Theory 
of Dynamical Systems, Cambridge University Press, 1995.
 
\bibitem[KO]{KO}Y. Katznelson, D. Ornstein, 
The differentiability of the conjugation of 
certain diffeomorphisms of the circle, Ergodic Theory Dynam. Systems 
9(1989), no.4, 643-680.

\bibitem[KS]{KS}K. M. Khanin, Ya. G. Sinai, New proof of M. Herman's 
theorem, Commun. Math. Phys. 112 (1987), 89-101.

\bibitem[M]{M}J. Milnor, Dynamics in one complex variable:
Introductory lectures, Math.ArXiv, math.DS/9201272. 
http://front.math.ucdavis.edu/math.DS/9201272

\bibitem[MS]{MS}W. de Melo, S. van Strien, One-dimensional Dynamics. 
Springer, New York, 1993.

\bibitem[PS]{PS}L. Polterovich, M. Sodin,  
A growth gap for diffeomorphisms of the interval, 
J.Anal.Math. 92(2004), 191-209.

\bibitem[St]{St}J. Stark, Smooth conjugacy and renormalization for 
diffeomorphisms of the circle, Nonlinearity 1 (4) (1988), 541-575.

\bibitem[Sw]{Sw}G. \'Swi\c atek, Rational rotation numbers for 
maps of the circle, Commun. Math. Phys. 119 (1988),109-128.

\bibitem[T]{T}F. Takens, Normal forms for certain singularities of vector 
fields, Ann. Inst. Fourier 23, no.2 (1973), 163-195.

\bibitem[W]{W}N. Watanabe, Growth sequences for flat diffeomorphisms of 
the interval, Nihonkai Math.J. 15(2004), no.2, 137-140. 

\bibitem[Y]{Y}J.-C. Yoccoz, Conjugaison diff$\acute{{\rm e}}$rentiable des 
diff$\acute{{\rm e}}$omorphismes du cercle dont le nomble de rotation 
v$\acute{{\rm e}}$rifie une condition diophantienne, 
Ann. Sci. \'Ecole Norm. Sup. 4 17 (1984), 333-359.  

\end{thebibliography}
\end{document}